# A Multi-server Scheduling Framework for Resource Allocation in Wireless Multi-carrier Networks

Ying Jun (Angela) Zhang, *Member, IEEE*
*The Department of Information Engineering, The Chinese University of Hong Kong*
*Shatin, New Territory, Hong Kong*
*Email: yjzhang@ie.cuhk.edu.hk*

Abstract- Multiuser resource allocation has recently been recognized as an effective methodology for enhancing the power and spectrum efficiency in OFDM (orthogonal frequency division multiplexing) systems. It is, however, not directly applicable to current packet-switched networks, because (1) most existing packet-scheduling schemes are based on a single-server model and do not serve multiple users at the same time; and (2) the conventional separate design of MAC (medium access control) packet scheduling and PHY (physical) resource allocation yields inefficient resource utilization. In this paper, we propose a cross-layer resource allocation algorithm based on a novel multi-server scheduling framework to achieve overall high system power efficiency in packet-switched OFDM networks. Our contribution is four fold. First, we propose and analyze a MPGPS (multi-server packetized general processor sharing) service discipline that serves multiple users at the same time and facilitates multiuser resource allocation. Second, we present a MPGPS-based joint MAC-PHY resource allocation scheme that incorporates packet scheduling, subcarrier allocation, and power allocation in an integrated framework. Third, by investigating the fundamental tradeoff between multiuser-diversity gain and queuing performance, we present an A-MPGPS (adaptive MPGPS) service discipline that strikes an optimal balance between power efficiency and queuing performance. Finally, we extend MPGPS to an O-MPGPS (opportunistic MPGPS) service discipline to further enhance the resource utilization efficiency. Through analysis, we prove that the proposed MPGPS, A-MPGPS, and O-MPGPS schemes can serve users in a way that approximates the ideal GPS (general processor sharing) service discipline, and hence guarantee QoS and fairness. Through simulations, we show that the MPGPS-, A-MPGPS-, and O-MPGPS-based cross-layer resource allocation algorithms significantly enhance system power efficiency compared to conventional resource allocation schemes.

Corresponding Author: Ying Jun (Angela) Zhang
Email: yjzhang@ie.cuhk.edu.hk
Tel:/Fax: (852) 2609 8465/2603 5032

This work is supported in part by the CUHK Direct Grant for Research #2050370.

# I. INTRODUCTION

The explosive growth in the deployment of wireless networks is expected to generate a tremendous increase in the demand for already scarce radio resources such as spectrum and energy. It is therefore essential to achieve efficient resource utilization in future broadband wireless networks. As one of the prime air interfaces for future broadband wireless networks, OFDM (orthogonal frequency division multiplexing) is well known for its flexibility of allowing single- and multi-user adaptive resource allocation to significantly enhance the efficiency of resource utilization [1-4]. In particular, OFDM multiuser resource allocation takes advantage of system diversities in the time, frequency, and user domains through dynamic allocation of subcarrier, rate, and power. Originally designed for cellular systems, most existing multiuser resource allocation schemes ignore the randomness in packet arrival and assume that the buffers are always infinitely backlogged [1-4]. As a result, they fail to appropriately allocate resources in wireless networks where traffic arrives in a random manner. For example, the PHY (physical) layer may allocate subcarriers and power to a user who has very little traffic in its queue. This causes a waste of resource. Moreover, the current algorithms usually assume a bit-by-bit scheduling, and hence are not directly applicable to packet-switched networks, where data is transmitted on a packet-by-packet basis.

On the other hand, packet scheduling schemes have been proposed for packet-switched networks to exploit the dynamics in traffic arrivals and ensure a fair and efficient bandwidth allocation [5-6]. Mainly confined to the MAC (medium access control) and network layers, these schemes are typically designed independently of the PHY layer, and hence are sometimes not compatible with advanced PHY-layer techniques such as OFDM multiuser resource allocation. For example, most scheduling schemes adopt a single-server model, which assumes that only one packet is served at a time [5-6]. However, to implement multiuser resource-allocation schemes such as those proposed in [1-4], multiple packets from different users should be served at the same time. It is therefore more appropriate to model the parallel subcarriers of OFDM systems as multiple servers, which allows simultaneous transmission of multiple packets. Unfortunately, in contrast to single-server scheduling, there is relatively little work on multi-server scheduling. Existing multi-server service disciplines usually assume that a packet is served by one server only and the busy periods of all servers do not coincide [7-8]. This assumption, however, does not apply to multiuser OFDM systems, where a packet can occupy more than one subcarrier and the busy periods of all subcarriers coincide.

Considering the sub-optimality of the separate design of packet scheduling and multiuser resource allocation, we propose a cross MAC-PHY layer resource allocation algorithm in this paper. The objective is to maximize system power efficiency while guaranteeing QoS (quality of service) and fairness. To capture the aforementioned features of multiuser OFDM systems, we propose and analyze a novel *multi-server* scheduling framework that assumes that (i) each packet can be served by a number of servers at the same time; (ii) multiple packets can be served simultaneously; and (iii) the busy periods of all servers coincide. Such a framework facilitates the exploitation of multiuser diversity in both MAC and PHY layers. Moreover, the PHY-layer power-and-subcarrier allocation is jointly optimized with the packet scheduling to achieve an overall good system performance. To the best of our knowledge, this is a first paper to propose and investigate a multi-server service discipline that is tailored for the multiuser resource allocation in OFDM networks.

Our contribution is four fold. First, we propose a novel MPGPS (multi-server packetized generalized processor sharing) service discipline. The performance of the proposed MPGPS service discipline is analyzed from the viewpoint of worst-case packet delay and throughput. The results show that MPGPS can serve packets in a way that approximates the ideal GPS (generalized processor sharing) service discipline, and hence ensures fairness and QoS in a system with random traffic arrivals. Second, built on the MPGPS service discipline, we formulate a joint MAC-PHY resource allocation problem that incorporates packet scheduling, subcarrier allocation, and power allocation in an integrated framework. By exploiting the special structure of the problem, we show that the optimal solution can be obtained in polynomial time. Our results will show that the power efficiency can be dramatically improved if resources are allocated according to the optimal solution. Meanwhile, fairness is guaranteed through both MPGPS scheduling in the MAC layer and the proportional subcarrier allocation in the PHY layer.

Based on the MPGPS resource allocation framework and given the inherent diversities in wireless networks, our next two contributions are to explore different dimensions to further improve the system power efficiency and service quality. Our third contribution is the investigation of a fundamental tradeoff between multiuser diversity and queuing performance, and the proposal of an A-MPGPS (adaptive MPGPS) based resource allocation algorithm that strikes an optimal balance between the two aspects. Our results show that the A-MPGPS algorithm effectively captures multiuser diversity without degrading the queuing performance unnecessarily. Our fourth contribution is the extension of MPGPS to an opportunistic service discipline, namely O-MPGPS (opportunistic MPGPS). The O-MPGPS scheme "opportunistically" swaps the order in which packets are transmitted to take advantage of the time variation of the channel. Most existing opportunistic scheduling schemes only provide services to users that temporarily experience good channels [9-11]. Hence, they fail to ensure short-term delay, throughput, and fairness to those who are stuck in deep fading channels. In contrast, the proposed O-MPGPS scheme guarantees QoS to each user regardless of the underlying channel conditions through (i) explicit control of the extent to which the scheduling is opportunistic and (ii) PHY-layer power control.

Through analysis, it is proved that the deviation of QoS relative to the non-opportunistic MPGPS is bounded. Our numerical results show that the O-MPGPS algorithm greatly outperforms the non-opportunistic MPGPS and A-MPGPS algorithms as well as the conventional resource allocation schemes.

The rest of this paper is organized as follows. Section II introduces the system model. In Section III, a novel MPGPS service discipline is proposed. Analysis of the worst-case delay and throughput is provided. The MPGPS-based cross layer resource allocation problem is then formulated in Section IV. Low-complexity algorithms that solve the joint optimization problem in polynomial time are presented. In Section V, we extend the proposed MPGPS to an A-MPGPS based resource allocation algorithm, which strikes an optimal balance between queuing performance and multiuser diversity. In Section VI, we extend MPGPS to an opportunistic service discipline, namely O-MPGPS, which further improves the system power efficiency while guaranteeing short-term QoS for each user. The performance of the proposed algorithms is evaluated by simulation in Section VII. Finally, the paper is concluded in Section VIII.

## II. SYSTEM MODEL

### *2.1 System Structure*

This paper considers the downlink transmission of an OFDM network with $K$ users and $N$ subcarriers. The transmitter structure at the BS (base station) is shown in Fig. 1. Packets arrive at the BS according to a random arrival process. Upon arrival, packets from different users are buffered in separate queues. When the channel becomes idle, the BS schedules packet transmission and allocates subcarriers and power to the selected packets according to the proposed resource allocation scheme. The packets are then modulated onto the assigned subcarriers and transmitted over the downlink channel using the allocated amount of power. If a packet is not correctly received, it stays in the queue as a HOL (head of line) packet and will be retransmitted at a future time, given that the deadline of the packet is not exceeded.

In this paper, we assume that packets are of equal length and contain $L$ bits each, and that fixed constellation size of $2^r$ is adopted on each subcarrier. Let $g_k$ denote the number of packets that are selected for transmission from user $k$'s buffer. Then, it takes

$$S \geq \left\lceil L\sum_{k=1}^{K} g_k \Big/ Nr \right\rceil \tag{1}$$

OFDM symbols to transmit the selected packets, where $\lceil \cdot \rceil$ denotes the smallest integer that is larger than or equal to the value of the input parameter. For simplicity, we assume that

$$S \geq \left\lceil L\sum_{k=1}^{K} g_k \Big/ Nr \right\rceil = L\sum_{k=1}^{K} g_k \Big/ Nr \quad (2)$$

in the rest of the paper, with the first equality hold if and only if the busy periods of all subcarriers coincde. In the following, we refer to a burst of $S$ OFDM symbols as one frame.

The proposed cross-layer resource allocation algorithm jointly optimizes packet scheduling, subcarrier allocation, and power distribution. In the rest of this section, we introduce the constraints on subcarrier and power allocation. Packet scheduling will be discussed in more detail in the subsequent sections.

*2.2 Power Allocation*

The objective of power allocation is to smooth the variance in the wireless channel and maintain a sufficiently low PER (packet error rate), so that the wireless channel appears to be static and almost error-free to the MAC layer. This greatly simplifies the design and analysis of MAC scheduling, as the techniques developed for wired networks can now be applied [12]. Let $H_{k,n,s}$ denote the channel coefficient on the $n^{\text{th}}$ subcarrier of user $k$ during the $s^{\text{th}}$ symbol of the current frame. Let $\varepsilon_p$ denote the target PER and $\gamma_k(\varepsilon_p)$ denote the required received SNR for user $k$ to achieve the target. Then, the transmit power is adapted according to the following equation:

$$\frac{\left|H_{k,n,s}\right|^2 p_{k,n,s}}{N_0 B} = \gamma_k(\varepsilon_p) \ \forall n,k \quad (3)$$

where $p_{k,n,s}$ denotes the transmit power subcarrier $n$ of symbol $s$ if it is assigned to user $k$, $N_0$ is the single-sided power spectral density of the white noise, and $B$ is the bandwidth of one subcarrier. If we assume that the channel remains constant during one frame, then $H_{k,n,s}$ and $p_{k,n,s}$ do not change over $s$. For brevity, we drop the subscript $s$ in the rest of the paper.

*2.3 Subcarrier Allocation*

From Eqn. (3), it can be seen that each subcarrier is associated with different transmit power for different users. Hence, the subcarriers need to be optimally allocated to minimize the total power consumption. Let $c_{k,n,s}$ be the subcarrier allocation indicator, with $c_{k,n,s} = 1$ if subcarrier $n$ of symbol $s$ is allocated to user $k$, and $c_{k,n,s} = 0$ otherwise. To avoid CCI (co-channel interference), a subcarrier is allocated to at most one user. That is,

$$\sum_{k=1}^{K} c_{k,n,s} = 1 \, \forall n, s\,. \quad (4)$$

Moreover, the number of subcarriers allocated to each user is determined by the number of packets transmitted from its queue during the current frame. That is,

$$\sum_{s=1}^{S}\sum_{n=1}^{N} c_{k,n,s} = \frac{g_k L}{r} \; \forall k \, . \tag{5}$$

## III. MULTI-SERVER PGPS

In this section, we describe and analyze our proposed MPGPS service discipline for multiuser adaptive OFDM systems, wherein *N* subcarriers are modeled as *N* servers. In particular, we extend the single-server PGPS (packetized generalized processor sharing) [6] to the multi-server case, assuming that (i) multiple packets can be served simultaneously by different servers; (ii) each packet can also be served by more than one server; and (iii) the service discipline is work conserving, which means that no servers can be idle as long as there are packets to transmit. In this case, the busy periods of all servers coincide.

Let *M* denote the number of packets that are served at the same time if there are enough backlogged packets at the time of scheduling. Single-server PGPS is a special case of the proposed MPGPS when $M = 1$. The MPGPS discipline works in the following way. When the servers become idle at time $\tau$, the system picks, among all the packets that are queued in the system at time $\tau$, the first *M* packets that would complete service in the corresponding GPS system if no additional packets were to arrive after time $\tau$. If the total number of backlogged packets is less than *M*, then all the backlogged packets are picked. The total number of packets that are scheduled at time $\tau$ is given by

$$\sum_{k=1}^{K} g_k = \min\left( M, \sum_{k=1}^{K} \hat{Q}_k(\tau) \right) \tag{6}$$

where $\hat{Q}_k(\tau)$ denotes the user-*k* backlog (in unit of packets) at time $\tau$ under MPGPS.

Similar to single-server PGPS, the MPGPS discipline uses the concept of virtual time to track the progress of GPS. When a packet arrives, it is stamped with its virtual time finishing time. Once the servers become idle, the $\min\left( M, \sum_{k=1}^{K} \hat{Q}_k(\tau) \right)$ packets with the smallest virtual time finishing time are selected for transmission. The details of the virtual time implementation of MPGPS are described in Appendix A. In the rest of this section, we focus on the analysis of the MPGPS service discipline. The following analysis is under the assumption that packet transmissions are error free. This assumption is reasonable, since the received signal strength is controlled by power allocation and we can always maintain a sufficiently low PER.

We first examine how much later packets may depart from the system under MPGPS relative to under GPS. Let us refer to the packets that are served at the same time as one *batch*. Let $f_h^k$ denote the $k^{\text{th}}$ packet in batch *h*, $\hat{d}_h^k$ denote the departure time of $f_h^k$ under MPGPS, $d_h^k$ denote the corresponding departure time under GPS, and $a_h^k$ denote the time when the packet arrives. Likewise, define $W_k(t_1, t_2)$ to be the amount of user *k* traffic (in unit of bits) served during the time period ($t_1$, $t_2$) under GPS and $\hat{W}_k(t_1, t_2)$ to be that under MPGPS.

***Theorem 1*** **(*Delay Bound Guarantee*):** For all packets, the difference between the departure times under MPGPS and GPS is bounded as follows:

$$\hat{d}_h^k - d_h^k \le \frac{(2M-1)L}{Nr} \qquad \forall h, k \tag{7}$$

*Proof:* Obviously, we have

$$\hat{d}_h = \hat{d}_h^1 = \hat{d}_h^2 = \cdots = \hat{d}_h^{M_h} \tag{8}$$

where $M_h$ denotes the actual number of packets that are served in batch *h*. Then, the starting time of batch *h*, denoted by $b_h$, is given by

$$b_h = \hat{d}_h - \frac{M_h L}{Nr} \tag{9}$$

because the discipline is work conserving. Without loss of generality, we assume that

$$d_h^1 \le d_h^2 \le \cdots \le d_h^{M_h}. \tag{10}$$

Let *j* be the largest integer that is smaller than *h* such that there exist *k'* and *k* that satisfy $d_j^{k'} > d_h^k$.

Thus
$$d_j^{k'} > d_h^k > d_i^{k''} \ \forall j < i < h, k''. \tag{11}$$

This implies that $f_j^{k'}, f_j^{k'+1}, \cdots, f_j^{M_j}$ are transmitted before packets $f_{j+1}^1, \cdots f_h^1, \cdots, f_h^k$ under MPGPS but after all these packets under GPS. This happens only when $f_{j+1}^1, \cdots f_h^1, \cdots, f_h^k$ arrive later than $b_j$.

That is,
$$\min\{a_{j+1}^1, \cdots, a_h^1, \cdots, a_h^k\} > \hat{d}_j - \frac{M_j L}{Nr}. \tag{12}$$

Since GPS is work conserving, it takes $\left(\sum_{l=1}^{h-j-1} M_{j+l} + k\right)L \Big/ Nr$ time slots to finish transmitting packets $f_{j+1}^1, \cdots f_h^1, \cdots, f_h^k$, if no packets were to arrive before the finishing time. Therefore, the departure time of the last packet $f_h^k$ under GPS satisfies

$$d_h^k \ge \min\{a_{j+1}^1, \cdots, a_h^1, \cdots, a_h^k\} + \frac{\left(\sum_{l=1}^{h-j-1} M_{j+l} + k\right)L}{Nr} > \hat{d}_j - \frac{M_j L}{Nr} + \frac{\left(\sum_{l=1}^{h-j-1} M_{j+l} + k\right)L}{Nr}, \tag{13}$$

while the departure time under MPGPS satisfies

$$\hat{d}_h^k = \hat{d}_j + \frac{\sum_{l=1}^{h-j} M_{j+l} L}{Nr}. \tag{14}$$

Combining (13) and (14), we have

$$\hat{d}_h^k \le d_h^k + \frac{M_j L}{Nr} + \frac{(M_h - k)L}{Nr} \; \forall h,k. \tag{15}$$

Since $M_j \le M$ and $M_h \le M$, it is straightforward from (15) that

$$\hat{d}_h^k \le d_h^k + \frac{(2M-1)L}{Nr} \qquad \forall h,k \tag{16}$$

and the theorem follows. □

Packets can depart later in MPGPS than in GPS due to the following reasons. First, multiple packets of one user may be concurrently in service at different servers of MPGPS. Second, the packets of the users with different priority weights are transmitted at the same rate if they are selected by the scheduler at the same time. Third, a packet that belongs to the first $M$ packets that would complete service under GPS may not have arrived at the time of scheduling. The first two types of late departure do not occur in the single-server case, but the last type also exists in single-server PGPS. In the next theorem, we show that even if no packet arrives too late to be served in time, packets can still depart later in MPGPS than in GPS (unlike single server PGPS).

***Theorem 2*****:** Even if all packets are served according to the increasing order of their departure times under GPS, they can still depart later in MPGPS than in GPS and the difference in the departure times is bounded by

$$\hat{d}_h^k - d_h^k \le \frac{(M-1)L}{Nr} \qquad \forall h,k. \tag{17}$$

*Proof*: Due to the fact that both GPS and MPGPS are work conserving, and that the packets are served in MPGPS according to the exact order of their departure times in GPS, the departure time of packet $f_h^k$ under GPS satisfies

$$d_h^k \ge b_h + \frac{kL}{Nr}, \tag{18}$$

with the equality holding if and only if $k = M_h$ and there is no more backlog at time $d_h^{M_h}$.

Substituting Eqn. (18) to (9), we have

$$\hat{d}_h^k \le d_h^k + \frac{(M_h - k)L}{Nr} \; \forall h,k \tag{19}$$

and hence $$\hat{d}_h^k - d_h^k \le \frac{(M-1)L}{Nr} \; \forall h,k$$ □

Next, we investigate the throughput guarantee of MPGPS.

***Theorem 3* (*Throughput Guarantee*):** For all times $t$ and session $k$,

$$W_k(0,t) - \hat{W}_k(0,t) \le (2M-1)L \, . \tag{20}$$

*Proof*: The theorem can be proved in a way similar to the proof of [6, Theorem 2] with needed modifications. Hence, we omit the proof here. For the convenience of the readers, we give the proof in Appendix B.

***Corollary 1*:** For all times $t$ and user $k$

$$\hat{Q}_k(t) - Q_k(t) \le 2M - 1 \, , \tag{21}$$

where $Q_k(t)$ and $\hat{Q}_k(t)$ denote the queue length (in unit of packets) of user $k$ at time $t$ under GPS and MPGPS, respectively.

***Comment:*** When $M$=1, the proposed MPGPS reduces to the single-server PGPS in [6]. In this case, Theorem 1–3 reduces to

$$\hat{d} - d \le \frac{L}{Nr} \, , \tag{22}$$

$$\hat{d} \le d \, , \tag{23}$$

and

$$W_k(0,t) - \hat{W}_k(0,t) \le L \, , \tag{24}$$

respectively. These results are consistent with the results in [6] when the rate of the PGPS server is equal to $Nr$.

It can be seen from the above derived results that the MPGPS discipline can provide users with similar service qualities as in a GPS system. It is well known that the worst-case packet delay and backlog of GPS are bounded if the traffic is regulated by a leaky bucket scheme [6]. Theorem 1 and Corollary 1 can be used to translate the bounds on GPS worst-case packet delay and backlog to the corresponding bounds on MPGPS. This implies that QoS and fairness requirements are ensured by the MPGPS scheme. On the other hand, the maximum deviation of the packet delay, throughput, and queue length with respect to GPS increases with $M$, which implies that the queuing performance degrades when the number of simultaneously served packets increases. This is due to the increase in the service time resulting from parallel packet transmission.

Note that so far we have assumed that the channel is error free, and hence no packet retransmission is considered. When packet retransmission is taken into consideration, queuing performance is also related to packet error rate, which in turn depends on system power efficiency. Usually, system power efficiency improves with the number of simultaneously served users, thanks to the multiuser diversity gain. Hence, there exists an optimal $M$ that strikes an optimal balance between power efficiency and queuing performance. It will be discussed in more detail in Section V.

## IV. MPGPS-BASED OPTIMAL RESOURCE ALLOCATION

In this section, we formulate the MPGPS-based optimal resource allocation into a LIP (linear integer programming) problem. By exploiting the structure the LIP problem, we show that it can typically be solved in polynomial time. The result of this section will be used as a basis for the development of two cross-layer extensions in Section V and Section VI.

Given that $g_k$ packets are selected from user $k$'s buffer by the scheduler at time $\tau$, the optimal subcarrier and power allocation that minimizes the power consumption per bit transmission is formulated into

$$\min_{c_{k,n,s},\, p_{k,n}} \sum_{s=1}^{S}\sum_{k=1}^{K}\sum_{n=1}^{N} \frac{p_{k,n} c_{k,n,s}}{L \sum_{k=1}^{K} g_k} \tag{25}$$

subject to constraints given by Eqn. (3)-(5) and

$$c_{k,n,s} = \{0,1\}\ \forall n,k,s, \tag{26}$$

where $S$ is given by Eqn. (2) and $g_k$ is given by the scheduling decision of MPGPS at time $\tau$. For the problem to be feasible, the RHS (right-hand-side) of (5) must be integral regardless of $g_k$. In other words, $L/r$ must be an integer. This is generally true for typical packet lengths and modulation schemes.

By substituting Eqn. (3) and (6) into (25), the optimization problem is converted to a LIP problem with an objective function given by

$$\min_{c_{k,n,s}} \sum_{s=1}^{S}\sum_{k=1}^{K}\sum_{n=1}^{N} \alpha_{k,n} c_{k,n,s} \tag{27}$$

and constraints described in Eqn. (4), (5) and (26). In (27)

$$\alpha_{k,n} = \frac{\gamma_k(\varepsilon_p) N_o B}{\left|H_{k,n}\right|^2 \min\left(M, \sum_{k=1}^{K} \hat{Q}_k(\tau)\right) L} \tag{28}$$

are constants that can be calculated at time $\tau$. The optimization problem can then be solved by LIP methods [13] and the solution yields the optimal subcarrier allocation, specified by $c_{k,n,s}$, and power allocation, specified by

$$p_{k,n} = \gamma_k(\varepsilon_p) N_o B \Big/ \left|H_{k,n}\right|^2 \quad \forall n,k, \tag{29}$$

given packet scheduling $g_k$.

The computational complexity of solving a LIP problem grows exponentially with the number of integer variables. In the above formulation, there are $NSK$ integer variables $c_{k,n,s}$, and hence the complexity is prohibitively high. Fortunately, the structure of the formulated problem allows us to significantly reduce the computational complexity.

*4.1 Problem-size reduction*

Note that the coefficients in the optimization problem remain constant for all *s*. Hence, we can divide the *S* OFDM symbols into small groups, and the optimal solution for one group applies to all the other groups. For example, when each group contains 1 OFDM symbol, the optimization problem for resource allocation in one group becomes

$$\min_{c_{k,n}} \sum_{k=1}^{K}\sum_{n=1}^{N} \alpha'_{k,n}\, c_{k,n} \tag{30.1}$$

Subject to:

$$\sum_{k=1}^{K} c_{k,n} = 1 \;\forall n \tag{30.2}$$

$$\sum_{n=1}^{N} c_{k,n} = \frac{g_k N}{\min\left( M, \sum_{k=1}^{K} \hat{Q}_k(\tau) \right)} \;\forall k \tag{30.3}$$

$$c_{k,n} = \{0,1\} \;\forall n,k \tag{30.4}$$

where

$$\alpha'_{k,n} = \frac{\gamma_k(\varepsilon_p) N_o B}{\left| H_{k,n}(\tau) \right|^2 Nr}. \tag{31}$$

Obviously, the size of the LIP problem reduces to *NK* integer variables. Note that the problem is feasible only when the RHS of (30.3) is integer. Otherwise, we have to increase the group size until the problem becomes feasible. In general, if the group size is *G*, the optimal resource allocation in one group is re-formulated into

$$\min_{c_{k,n,s},\, p_{k,n}} \sum_{s=1}^{G}\sum_{k=1}^{K}\sum_{n=1}^{N} \alpha''_{k,n} c_{k,n,s} \tag{32.1}$$

$$\sum_{k=1}^{K} c_{k,n,s} = 1 \;\forall n, s = 1 \ldots G \tag{32.2}$$

$$\sum_{s=1}^{G}\sum_{n=1}^{N} c_{k,n,s} = \frac{g_k GN}{\min\left( M, \sum_{k=1}^{K} \hat{Q}_k(\tau) \right)} \;\forall k \tag{32.3}$$

$$c_{k,n,s} = \{0,1\} \;\forall n,k, s = 1 \ldots G \tag{32.4}$$

with

$$a''_{k,n} = \frac{\gamma_k(\varepsilon_p) N_o B}{\left| H_{k,n}(\tau) \right|^2 GNr} \tag{33}$$

To achieve feasibility while making the problem size as small as possible, *G* is set to be the smallest integer that divides *S* while ensuring that the RHS of (32.3) is integer.

*4.2 Integer constraint relaxation*

Even if $G$=1, the computational complexity of solving the LIP problem still grows exponentially with $NK$, where $N$ is typically at the order of $10^2$~$10^3$. In this subsection, we show that even if we ignore the integer constraints on $c_{k,n,s}$, the solutions are still guaranteed to be integers. As a result, the LIP problem can be converted to a LP (linear programming) problem by relaxing the integer constraints. The computational complexity is therefore drastically decreased.

***Proposition 1:*** The optimal solution to problem (32.1-32.4) is guaranteed to be integers, even if the integer constraints on $c_{k,n,s}$ are relaxed and Constraint (32.4) is replaced by $0 \le c_{k,n,s} \le 1$.

*Proof:* The constraints (32.2) and (32.3) can be written into a matrix form

$$\mathbf{A}\mathbf{c} = \mathbf{b} \tag{34}$$

where $\mathbf{c} = \left[c_{1,1,1} \cdots c_{K,N,G}\right]^T$, $\mathbf{b} = \left[1,\cdots 1, \frac{g_1 GN}{\min\left(M, \sum_{k=1}^{K} \hat{Q}_k(\tau)\right)}, \cdots, \frac{g_K GN}{\min\left(M, \sum_{k=1}^{K} \hat{Q}_k(\tau)\right)}\right]^T$, and $\mathbf{A}$ contains the coefficients of the linear equations. Obviously, the determinant of each square submatrix of $\mathbf{A}$ is either 0, +1, or −1, which means that $\mathbf{A}$ is a unimodular matrix. Meanwhile, each element of $\mathbf{b}$ is an integer. In [13], it has been proved that the optimal solution of a LP problem is guaranteed to be integral if the constraint matrix $\mathbf{A}$ is unimodular and all elements of vector $\mathbf{b}$ are integers. A detailed proof can also be found in [12, Appendix A].

Now the optimization problem is simplified to a LP problem. Empirical experiments show that such a problem can typically be solved in polynomial time.

## V. ADAPTIVE MPGPS BASED RESOURCE ALLOCATION

The analysis in Section III indicates that the queuing performance degrades as the number of simultaneously served packets increases, assuming there is no packet retransmission. When packet retransmission is taken into account, the queuing performance is also affected by packet error rate that improves with system power efficiency. It is well known that a large number of simultaneously served users typically lead to high power efficiency in a system with multiuser resource allocation. As a result, increasing $M$ (the number of concurrently served packets) has two contrary effects. On the one hand, it degrades the queuing performance due to parallel transmission. On the other hand, it improves the queuing performance through improving power efficiency, if the simultaneously transmitted packets belong to different users. In view of this tradeoff, we propose in this section an A-MPGPS scheme, which adapts the number of concurrently served packets to strike an optimal balance.

The idea of A-MPGPS is as follows. Given a particular channel realization, serving more packets at a time may not lead to higher power efficiency, especially when the additional packet belongs to a user that perceives a deep fading channel. In this case, it is better to reduce the number of concurrently served packets to avoid unnecessary degradation of queuing performance. Therefore, the proposed algorithm adapts the value of *M* so that the deep-faded users are excluded from service, the multiuser diversity is effectively exploited, and the QoS is not unnecessarily degraded.

In the following, we describe the proposed A-MPGPS based resource allocation algorithm. In particular, at every scheduling time, we increase *M* gradually until increasing *M* does not improve the power efficiency any more.

***A-MPGPS Based Energy-Efficient Resource Allocation***

*Initialization*:

Define $M_{\max}$ to be the maximal allowable value of *M*. First, let *M*=1 and pick, from the backlogged buffers, the packet with the smallest virtual time finishing time. Since there is now only one active user, optimal subcarrier-and-power allocation reduces to single user water-filling. Let $P^{\min}$ denote the resultant power consumption per bit. Likewise, let $g_k^{\min}$ denote the scheduling result.

*Step 1*: *Scheduling*

Set *M*=*M*+1. Select *M* packets with the smallest virtual time finishing time from the backlogged buffers and let $g_k$ denote the number of packets that are selected from user *k*'s buffer.

*Step 2*: *Subcarrier-and-power allocation*

Given $g_k$, find the optimal solution to (25). Let $c_{k,n,s}^*$, $p_{k,n}^*$, and $P^*$ denote the corresponding optimal values of $c_{k,n,s}$, $p_{k,n}$, and *P*, respectively.

If $P^* < P^{\min}$, then let $c_{k,n,s}^{\min} = c_{k,n,s}^*$, $p_{k,n}^{\min} = p_{k,n}^*$, $P^{\min} = P^*$, and $g_k^{\min} = g_k$.

If $P^* \geq P^{\min}$ or $M = M_{\max}$, then go to Step 3. Otherwise, go back to Step 1.

*Step 3*: *Termination*

The resource allocation solution is given by $c_{k,n,s}^{\min}$, $p_{k,n}^{\min}$, and $g_k^{\min}$.

## VI. Opportunistic MPGPS Based Resource Allocation

So far, the order in which packets are served has strictly followed the increasing order of the virtual time finishing time of packets under GPS if no packets would arrive too late to be scheduled in time. In this section, we extend MPGPS to an O-MPGPS service discipline that exploits the time-domain diversity by "opportunistically" swapping the order in which packets are transmitted. By jointly designing O-MPGPS scheduling and subcarrier-and-power allocation, the proposed system is able to make full use of system diversities in the time, frequency, and user domains.

One disadvantage of existing channel-state-dependent scheduling algorithms is that a user may lose excessive share of service if it happens to be stuck in a deep fading channel for a long time. When this happens, the lag in the user's queue keeps increasing until packets are dropped due to delay-bound violation. In contrast, the sum of the lags (as well as leads) in all queues is bounded in the proposed O-MPGPS system. As a result, the share of service to each user is guaranteed. This is achieved by limiting the extent to which O-MPGPS is opportunistic. In other words, instead of choosing the best $M$ packets out of all backlogged packets, we choose from a *subset* of the backlogs with the smallest virtual time finishing time.

The algorithm works as follows. When the servers become idle at time $\tau$, the system identifies the first $U$ packets that would complete the service in the corresponding GPS system if no additional packets were to arrive after time $\tau$. Out of the $U$ packets, the $M$ packets that maximize the power efficiency are selected. The selection of the $M$ packets (or the decision of $g_k$) should be jointly optimized with subcarrier allocation and power allocation. If out of the $U$ packets $U_k$ belongs to user $k$'s queue, then $g_k$ is subject to the constraint $g_k \in \{0,1,\cdots,U_k\}$. The optimization problem can then be formulated as follows, with $c_{k,n,s}$ and $g_k$ being the variables.

$$\min_{c_{k,n,s},g_k} \sum_{s=1}^{S}\sum_{k=1}^{K}\sum_{n=1}^{N}\alpha_{k,n}c_{k,n,s} \tag{35}$$

subject to constraints given by (4), (5), (26), and

$$\sum_{k=1}^{K} g_k = \min\left(M, \sum_{k=1}^{K}\hat{Q}_k(\tau)\right) \tag{36}$$

$$g_k \in \{0,1,\cdots,U_k\} \forall k \tag{37}$$

where $S$ and $\alpha_{k,n}$ is given by (2) and (28), respectively.

In the following, we derive the upper bound on the aggregate lag (or lead) of the O-MPGPS system with respect to MPGPS. Note that the aggregate lag is equal to the aggregate lead, since the service discipline is work-conserving. Assume that among the $U$ candidate packets, there are $G_{lag}$ packets that are to be transmitted before the current time under MPGPS (referred to as lagging packets), $G_{sync}$ packets that are to be transmitted at the current time under MPGPS (referred to as in-sync packets), and $G_{lead}$ packets that are to be transmitted after the current time under MPGPS (referred to as leading packets). Also assume that the O-MPGPS scheme picks $M_1$ out of $G_{lag}$, $M_2$ out of $G_{sync}$, and $M_3$ out of $G_{lead}$ packets to transmit. Obviously, we have

$$G_{lag} + G_{sync} + G_{lead} = U \tag{38}$$

and $$M_1 + M_2 + M_3 = M \, . \tag{39}$$

Note that the $G_{sync} - M_2$ un-served in-sync packets will add to the aggregate lag at the next scheduling time. As a result, the net increase in the aggregate lag, denoted by *A*, is given by

$$A = G_{sync} - M_2 - M_1 \, . \tag{40}$$

Since the lagging packets have the smallest virtual time finishing time, $G_{lag}$ is equal to the aggregate lag as long as the aggregate lag is no larger than *U*. In fact, as we will prove in Theorem 4, the aggregate lag (or lead) is always no larger than $U - M$ in the O-MPGPS system.

***Theorem 4:*** The aggregate lag (or lead) of the O-MPGPS scheme (relative to MPGPS) is upper bounded by $U - M$ .

*Proof:* To prove Theorem 4, we show that when the aggregate lag is equal to $U - M - \delta$, the net increase *A* is always no larger than $\delta$ . Assume that $G_{lag} = U - M - \delta$ . Then, we have

$$G_{sync} + G_{lead} = M + \delta \tag{41}$$

and $$M_1 = M - M_2 - M_3 \tag{42}$$

Substituting (41) to (42), we have

$$M_1 = G_{sync} - M_2 + G_{lead} - M_3 - \delta \geq G_{sync} - M_2 - \delta \, . \tag{43}$$

The last inequality holds due to the fact that $M_3 \leq G_{lead}$ . Hence, we have

$$A = G_{sync} - M_2 - M_1 \leq \delta \, . \tag{44}$$

□

***Comment:*** In O-MPGPS, a user is guaranteed to receive service when its loss of share of bandwidth (i.e., the number of lagging packets) reaches an upper bound. The signal reception quality is guaranteed through PHY-layer power allocation.

***Corollary 2:*** Since the aggregate lag is upper-bounded in O-MPGPS, the deviation of QoS and fairness with respect to MPGPS is also upper-bounded. Interested readers are referred to [12] for the proof of a similar statement.

## VII. NUMERICAL RESULTS

In this section, the performance of the proposed MPGPS, A-MPGPS, and O-MPGPS based resource allocation algorithms is investigated by simulations. In particular, we first examine how system performance is affected by *M*, the number of concurrently served packets, and *U*, the opportunistic freedom. Then, we demonstrate the significant performance improvement of the proposed scheme compared with the conventional systems without joint MAC-PHY optimization. Throughout this section, we assume that the number of users *K*=10, the number of subcarriers *N*=64, the packet length $L$=1024, the OFDM symbol duration is 200 $\mu s$, and QPSK (quadrature phase shift keying) modulation is adopted on each subcarrier (i.e., $r$=2). For simplicity, all users have the same priority, though the proposed algorithms and analyses are applicable to any other priority settings. Users are randomly located within a circular cell with a radius of 50 meters. We assume a path loss exponent of 4 and a log-normal shadowing with standard deviation of 6dB. The channel is Rayleigh faded with an exponential power delay profile. For each simulation, we ran 50 replications, with each lasting for 28 hours.

*A.* Effect of *M* and *U* on system performance

In this subsection, we investigate the effect of *M* and *U* on power efficiency, packet delay, packet loss rate, and fairness.

In Fig. 2, we demonstrate the significant power-efficiency enhancement brought about by the proposed schemes when the users are infinitely backlogged and the target BER is $10^{-6}$ (which corresponds to $\gamma_k = 13.5401 dB$). For comparison, the average power consumption of conventional OFDM systems with a single-server PGPS service discipline (i.e., $M$=1 and $U$=1) is also plotted. It can be seen that the transmit power needed to ensure the signal strength at the receiver is significantly reduced by more than 10dB in the proposed schemes when $M$ (for MPGPS and O-MPGPS) or $M_{\max}$ (for A-MPGPS) is larger than one. This is because the simultaneous packet transmission facilitates the exploitation of multiuser diversity. By allowing opportunistic packet scheduling, the O-MPGPS scheme outperforms MPGPS and A-MPGPS by another 2dB (when $U = 10$). This is because the order of packet transmission can now be swapped to take advantage of the channel variation in the time domain. An interesting observation is that the A-MPGPS scheme, where the number of concurrently served packets is typically less than $M_{\max}$, is able to achieve a similar multiuser diversity gain as the MPGPS scheme with $M = M_{\max}$.

In the following, we investigate a practical scenario where the traffic load is finite. Assume that the source of each flow generates data packets according to a Poisson process with an average rate of 63*k*bps. A packet is dropped if it is not successfully transmitted within $40\,ms$ after its arrival.

First, we investigate the impact of $M$ on the average packet delay in Fig. 3. It is not surprising to see that the average packet delays of MPGPS, A-MPGPS, and O-MPGPS increase with $M$ (or $M_{\max}$), as proved in Section III. An interesting observation is that the A-MPGPS scheme effectively reduces the average packet delay compared to MPGPS by avoiding unnecessary parallel packet transmissions. There is more than 13% reduction when $M$ (or $M_{\max}$) is equal to 6, for instance. This observation, together with that in Fig. 2, shows that the A-MPGPS scheme can achieve both lower average packet delay and slightly higher power efficiency compared with MPGPS.

A close observation of Fig. 3 indicates that opportunistic transmission (i.e., letting $U > M$ ) hardly affects the mean system performance such as average packet delay. Instead, as we will show in the following figures, $U > M$ has an impact on fairness as well the tail probability of system performance, such as packet loss rate due to delay bound violation. In Fig. 4, we investigate fairness in terms of the maximum difference between the normalized services received by any two users *i* and *j* during any time interval $(t_1, t_2)$ in which both users are continuously backlogged. This is given by

$$\max\left|\frac{W_i(t_1,t_2)}{\phi_i} - \frac{W_j(t_1,t_2)}{\phi_j}\right|,$$

where $\phi_i$ and $\phi_j$ denote the priority weight of user *i* and *j*, respectively. In our simulations, $\phi_i = \phi_j = 1$. From the figure, we can see that fairness degrades with the increase of *U*. This is because the original order of packet transmission under MPGPS, which preserves fairness, is changed according to the opportunistic service discipline. The larger the *U*, the more freedom there is to swap the order of packet transmission. Meanwhile, it can be seen that given a fixed *U*, fairness improves with the increase of *M*. This is because a large *M* decreases the freedom of opportunistic scheduling. In particular, when *M* approaches *U*, the service discipline reduces to a non-opportunistic one. Similar conclusions can be drawn from Fig. 5, where packet loss rate is plotted as a function of *U* and *M*.

The power efficiency of O-MPGPS is affected by both *U* and *M*, as shown in Fig. 6. When *M* is fixed, it is not surprising to see that the power efficiency improves with the increase of *U*. However, for a fixed *U*, the effect of *M* is more complicated. On the one hand, a large *M* leads to a high multiuser diversity gain resulting from dynamic subcarrier-and-power allocation. On the other hand, it results in reduced freedom of opportunistic scheduling. Hence, there is a tradeoff between the multiuser diversity gain and opportunistic scheduling gain. It can be seen from Fig. 6 that when $U \le 4$, the average transmit power decreases as *M* increases. When *U* is equal to 5 or 6, the power consumption is minimized at $M = 4$ under the given traffic load and channel model.

The above figures have shown that *U* and *M* have an impact on both power efficiency and queuing performance of O-MPGPS. To achieve an overall good performance, it is necessary to set the parameters appropriately so that they strike an optimal balance. In Fig. 7, we show that for a fixed *U*, there exists an optimal *M* that outperforms other *M*'s in both power efficiency and queuing performance. In the figure, packet loss rate is plotted against the average transmit power when $U = 6$ and $M$ varies from 1 to 6. It can be seen that $M = 5$ achieves the pareto optimum. That is, it yields both higher power efficiency and lower packet loss rate compared to other values of $M$.

*B.* Performance comparison with conventional schemes

Having observed the impact of $M$ and $U$ in the proposed service disciplines, we proceed to demonstrate the advantage of the proposed schemes over (i) OFDM systems with a conventional single-server PGPS service discipline and (ii) traditional subcarrier-and-power allocation schemes that assume infinite backlog at each queue and do not apply statistical multiplexing. In particular, (i) corresponds to a special case of the proposed multi-server PGPS schemes when $M = 1$ and $U = 1$. In (ii), subcarriers are always distributed among all users, even though some users may not have enough packets to transmit during certain timeslots.

In Fig. 8, throughput, which is defined as the number of successfully transmitted packets within one OFDM symbol, is plotted against the average transmit $E_b/N_o$ (per-bit energy normalized by the power spectrum density of noise). The figure shows that when the transmit power is low, the service rate of the channel is also low due to the high packet error rate. In this case, the throughput is limited by the service rate of the channel. While in the high transmit-power region where packet error rate is low, the service rate is higher than the packet arrival rate. In this case, the throughput is determined by the amount of input traffic. Note that the system is stable only when the service rate of the channel is higher than the packet arrival rate.

In Fig. 9, packet loss rate is plotted against the average transmit power. When the system becomes unstable in the low transmit-power region, packet loss rate due to delay bound violation degrades drastically. From the figure, it can be seen that by exploiting the randomness in packet arrivals, both single-server PGPS and the proposed multi-server PGPS schemes yield a much lower packet loss rate compared to the conventional resource allocation schemes without statistical multiplexing. A dramatic decrease of more than 80% is observed in the figure.

Both Fig. 8 and Fig. 9 show that the conventional system with a single-server PGPS service discipline (which corresponds to $U = 1$ and $M = 1$) requires at least 22dB transmit $E_b/N_o$ to achieve a stable performance. By increasing $U$ from 1 to 2, we achieve a 2dB improvement in the power efficiency, thanks to the freedom of opportunistic scheduling. By introducing multi-server scheduling, the system performance is greatly improved due to the multiuser diversity gain. By serving 2 packets at the same time (A-MPGPS with $M_{\max}$=2), the power efficiency is dramatically improved by around 8dB compared to single-server PGPS. When $M_{\max}$ is increased to 6, another more than 1dB power gain is achieved. O-MPGPS achieves higher power efficiency than A-MPGPS, thanks to its capability of making full use of system diversities in the time, frequency, and user domains. With $U = 6$ and $M = 5$ (which is the optimal $M$ according to Fig. 7), O-MPGPS greatly outperforms single-server PGPS and A-MPGPS with $M_{max}$=6 by around 13dB and 3dB, respectively.

In this section, we have demonstrated the advantage of the proposed MPGPS, A-MPGPS, and O-MPGPS schemes over conventional packet scheduling and resource allocation schemes. We have also investigated the tradeoff between power efficiency and QoS such as packet delay, packet loss rate, and fairness. In practice, appropriate $U$ and $M$ can be selected according to the QoS requirements of the application, the load of the system, the power budget, and the affordable computational cost.

## VIII. Conclusions

In this paper, we have (i) proposed and analyzed a MPGPS service discipline that facilitates the exploitation of multiuser diversity in packet-switched OFDM networks; (ii) extended MPGPS to an A-MPGPS scheme that strikes an optimal balance between the multiuser diversity gain and queuing performance; (iii) proposed and analyzed an O-MPGPS service discipline that further enhances system power efficiency by exploiting the time-domain diversity; and (iv) presented MPGPS-, A-MPGPS-, and O-MPGPS-based cross MAC-PHY layer resource allocation algorithms, which incorporate packet scheduling, subcarrier allocation, and power allocation into one integrated framework.

With regard to (i), we have proved that the MPGPS discipline serves the users in a way that approximates the ideal GPS service discipline. The maximum deviations of packet delay, throughput, and queue length with respect to GPS are upper-bounded. With regard to (ii), we tackled the fundamental tradeoff between the multiuser diversity gain and queuing performance based on the analysis in (i). Numerical results have shown that the proposed A-MPGPS scheme effectively captures the multiuser diversity gain, while noticeably reducing the average packet delay. With regard to (iii), we proved that the aggregate lead or lag with respect to non-opportunistic MPGPS, and hence the degradation of QoS, is upper bounded in the proposed O-MPGPS system. With regard to (iv), we exploited the special structure of the LIP problem and proposed reduce-complexity schemes that solve the problem in polynomial time. Numerical results have shown that the power efficiency is dramatically improved by the proposed joint MAC-PHY layer resource allocation. Compared with the conventional resource allocation that assume infinite backlog, the packet loss rate due to delay bound violation is drastically reduced by around 80%. Compared with single-server PGPS based resource allocation schemes, the A-MPGPS based resource allocation greatly improves power efficiency by 8dB and 9dB when $M_{max}$ is 2 and 6, respectively, while the power efficiency of the O-MPGPS based system is significantly enhanced by 13dB when $U = 6$ and $M = 5$.

## APPENDIX A: VIRTUAL TIME IMPLEMENTATION OF MPGPS

Similar to [6], we use the concept of Virtual Time to track the progress of GPS. The difference is that the service rate is now assumed to be *Nr* instead of 1. Let $t_1 = 0$ denote the time of the first arrival of a busy period and $t_j$ denote the time at which the $j^{th}$ event (either arrival or departure) happens. Likewise, let $B_j$ denote the set of backlogged users during the interval $(t_{j-1,}t_j)$, $V(t)$ denote the virtual time of the system, and $\phi_i$ denote the priority weight of user *i*. Consider any busy period, then $V(t)$ evolves as follows:

$$V(0) = 0 \tag{A1}$$

$$V(t_{j-1}+\tau)=V(t_{j-1})+\frac{Nr\cdot\tau}{\sum_{i\in B_j}\phi_i}\text{, for }\tau\le t_j-t_{j-1},\ j=2,3,\cdots \tag{A2}$$

When the $k^{\text{th}}$ user $i$ packet arrives at time $a_i^k$, its virtual time starting time $S_i^k$ and virtual time finishing time $F_i^k$ are calculated by

$$S_i^k=\max\left(F_i^{k-1},V(a_i^k)\right)\text{ and }F_i^k=S_i^k+\frac{L}{\phi_i}. \tag{A3}$$

It is obvious from (A2) that the next virtual time update after $t$, if there are no more arrivals after time $t$, is give by

$$Next(t)=t+\left(F_{\min}-V(t)\right)\frac{\sum_{i\in B_j}\phi_i}{Nr}. \tag{A4}$$

APPENDIX B: PROOF OF THEOREM 3

Since the slopes of both $W_k$ and $\hat{W}_k$ vary between 0 and $Nr$, the difference $W_k(0,t)-\hat{W}_k(0,t)$ reaches its maximal value when user $k$ packets begin transmission under MPGPS. Let $t=b_j$ be the starting time of batch $j$ during which user $k$ receives service, $g_k$ be the number of user $k$ packets that are transmitted in batch $j$, and $\tau$ be the time at which all the $g_k$ packets complete service under GPS. Since packets of one user are served in a FIFO order under both schemes,

$$W_k(0,\tau)=\hat{W}_k\left(0,t+\frac{M_jL}{Nr}\right). \tag{A5}$$

Substituting $d_j^k=\tau$ and $\hat{d}_j^k=t+\frac{M_jL}{Nr}$ in (15), we have

$$\tau\ge t+\frac{M_jL}{Nr}-\frac{ML}{Nr}-\frac{(M_j-1)L}{Nr}=t-\frac{(M-1)L}{Nr}, \tag{A6}$$

which leads to

$$W_k\left(0,t-\frac{(M-1)L}{Nr}\right)\le W_k(0,\tau)=\hat{W}_k\left(0,t+\frac{M_jL}{Nr}\right). \tag{A7}$$

Since the slope of $W_k$ is at most $Nr$,

$$\begin{aligned}W_k(0,t)-(M-1)L&\le W_k\left(0,t-\frac{(M-1)L}{Nr}\right)\le\hat{W}_k\left(0,t+\frac{M_jL}{Nr}\right)=\hat{W}_k(0,t)+M_jL\\&\le\hat{W}_k(0,t)+ML\end{aligned} \tag{A8}$$

and hence

$$W_k(0,t)-\hat{W}_k(0,t)\le(2M-1)L. \tag{A9}$$

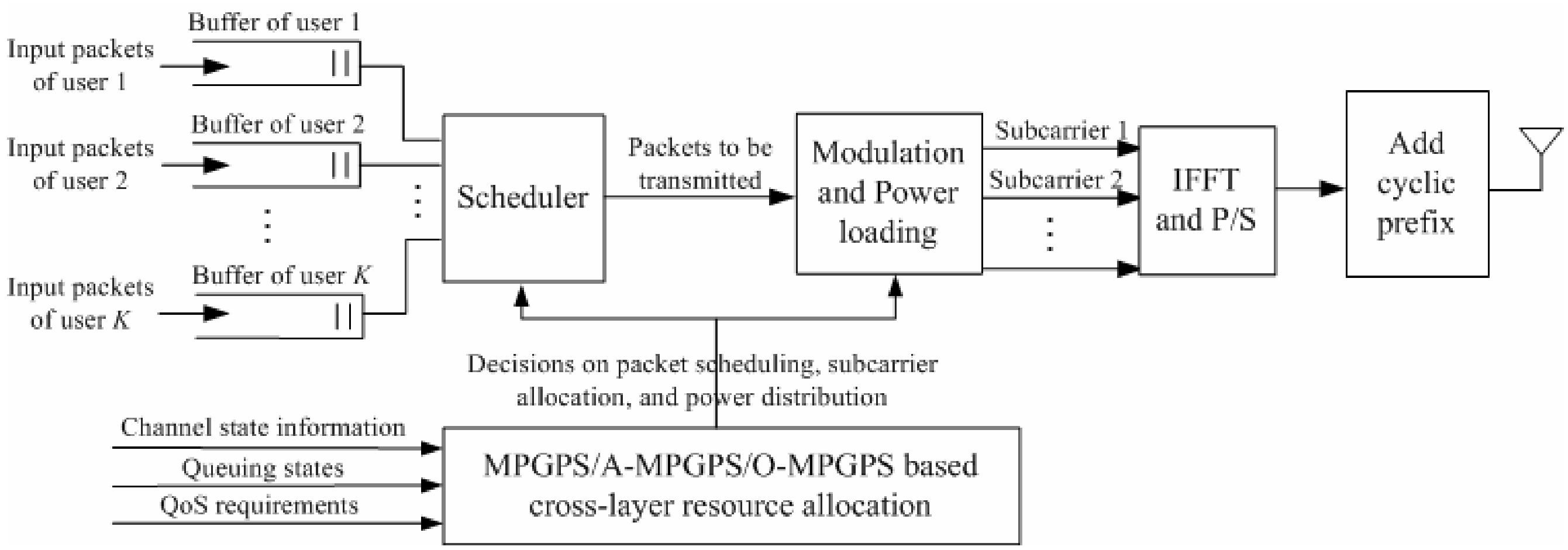

**Fig. 1:** Transmitter structure of a system with cross-layer MPGPS/A-MPGPS/O-MPGPS based resource allocation

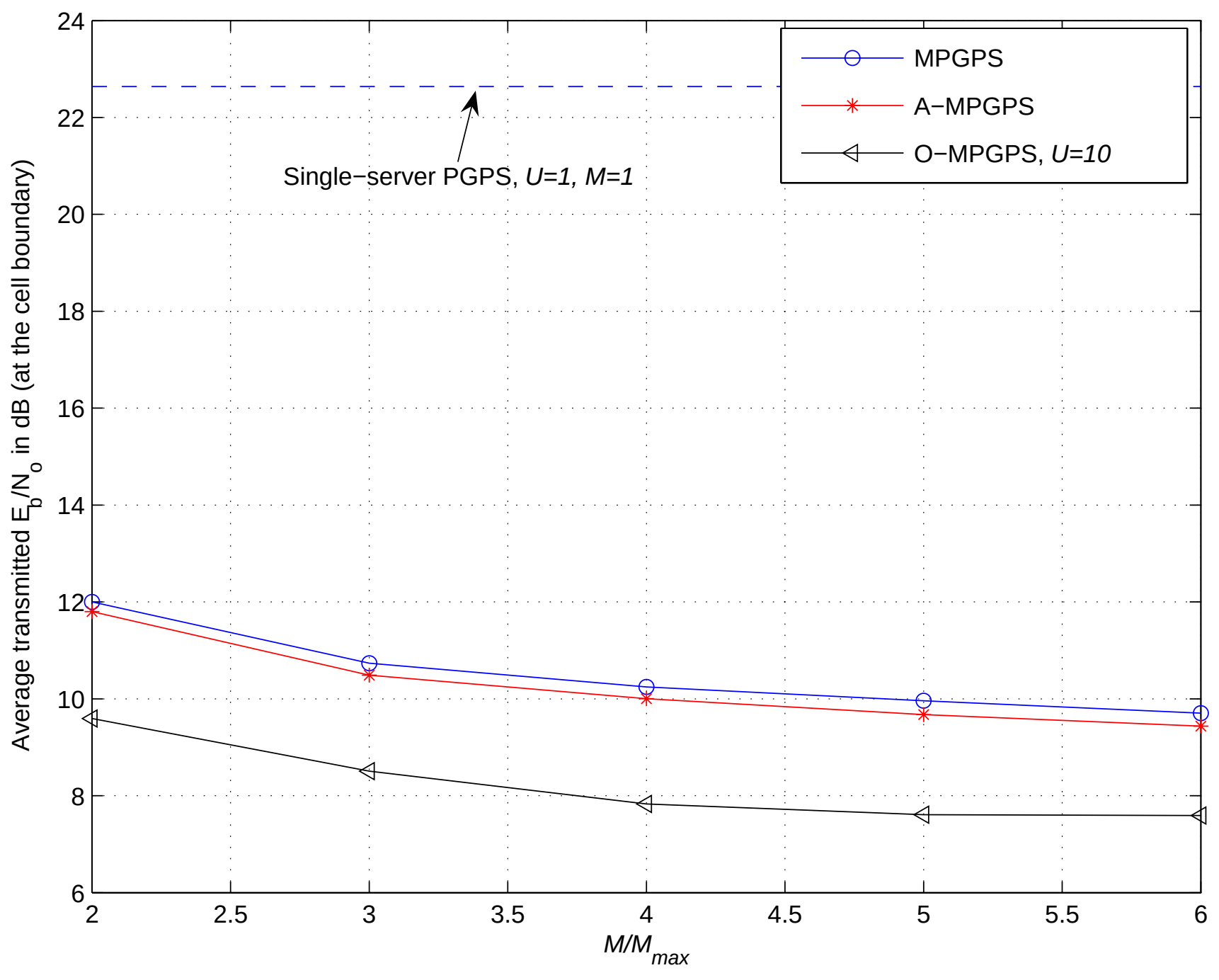

**Fig. 2:** Comparison of power efficiency between MPGPS, A-MPGPS, and O-MPGPS

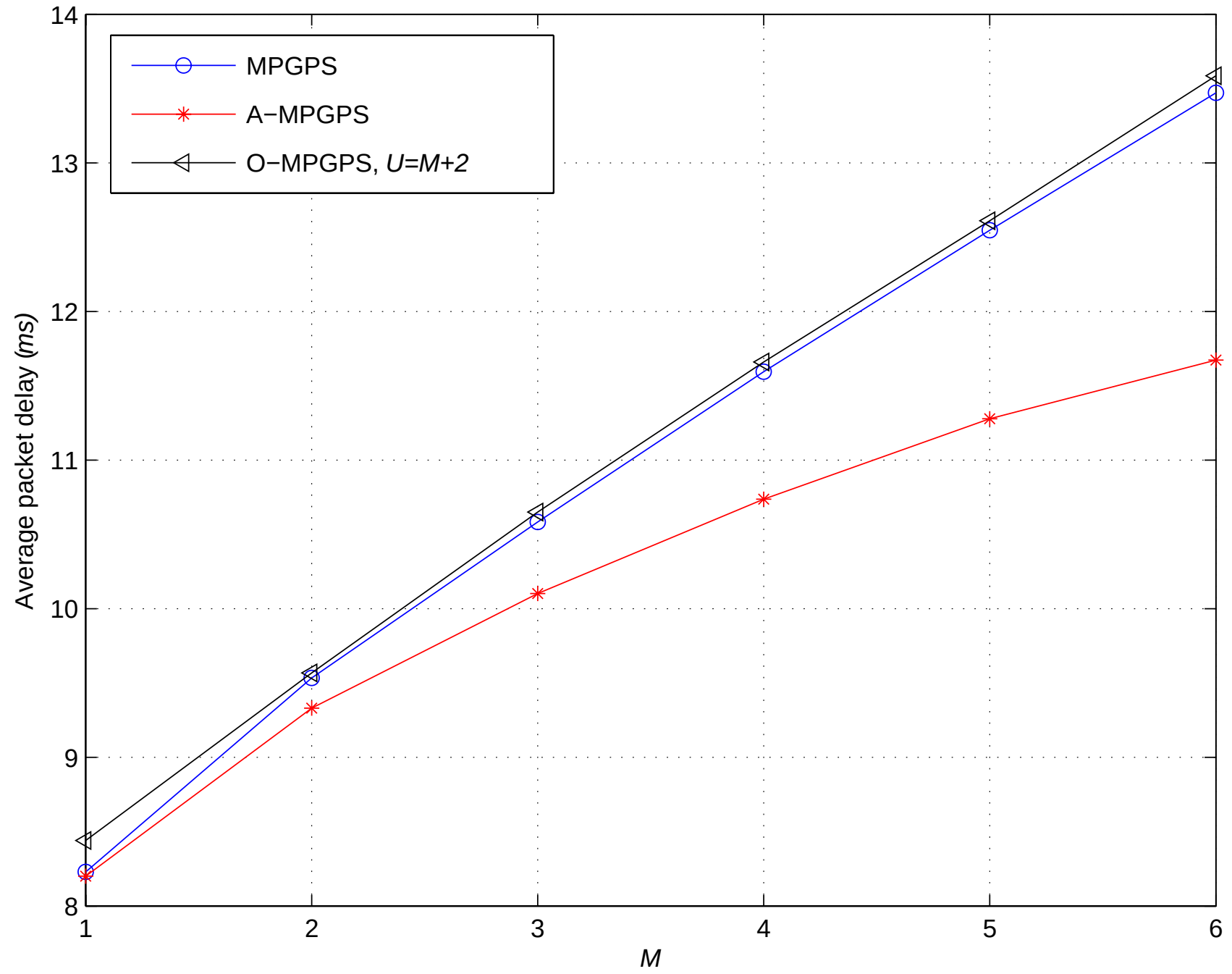

**Fig. 3:** Comparison of average packet delay between MPGPS, A-MPGPS, and O-MPGPS

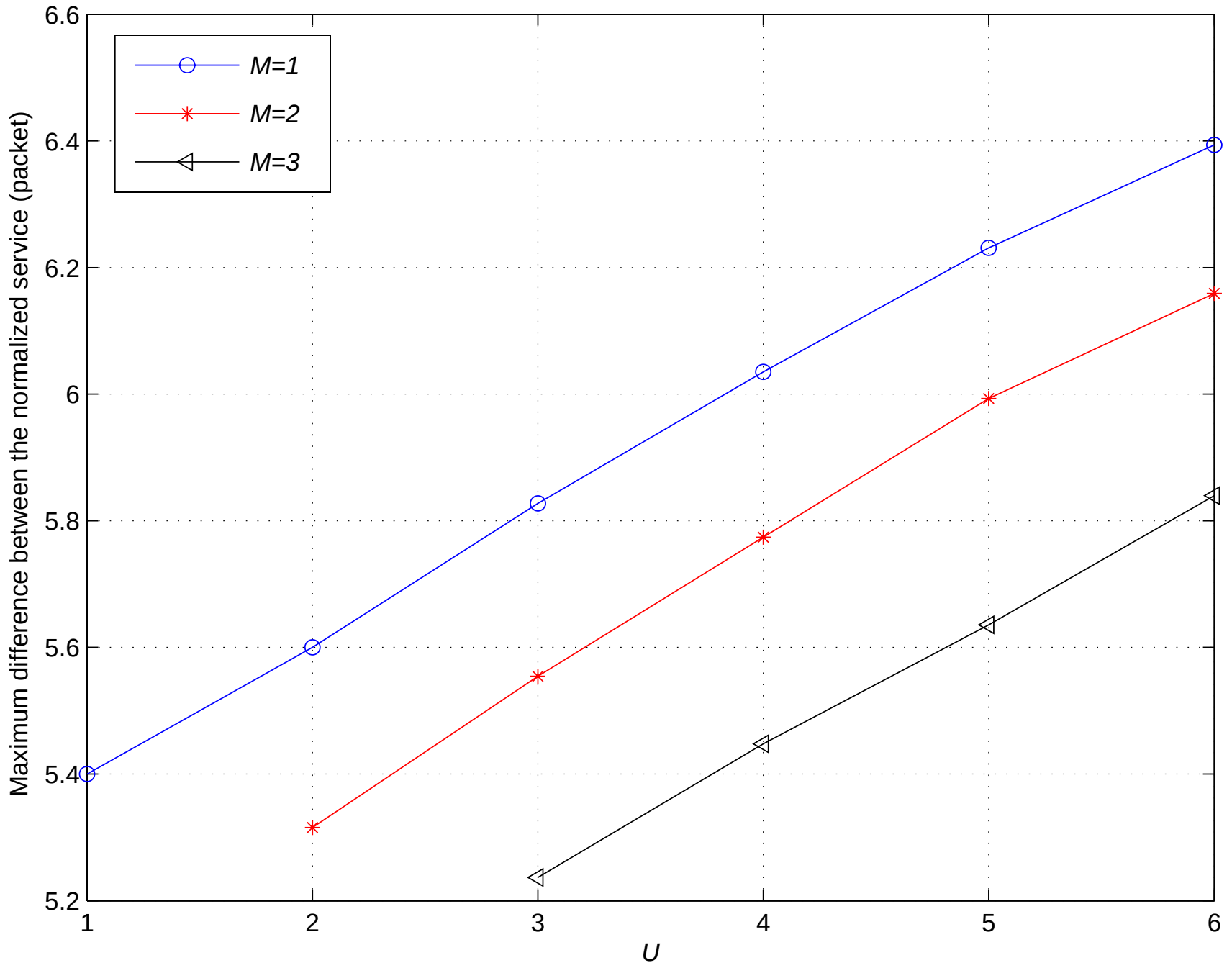

**Fig. 4:** Impact of *U* on fairness. $t_2 - t_1 = 100ms$, $\phi_i = \phi_j = 1$.

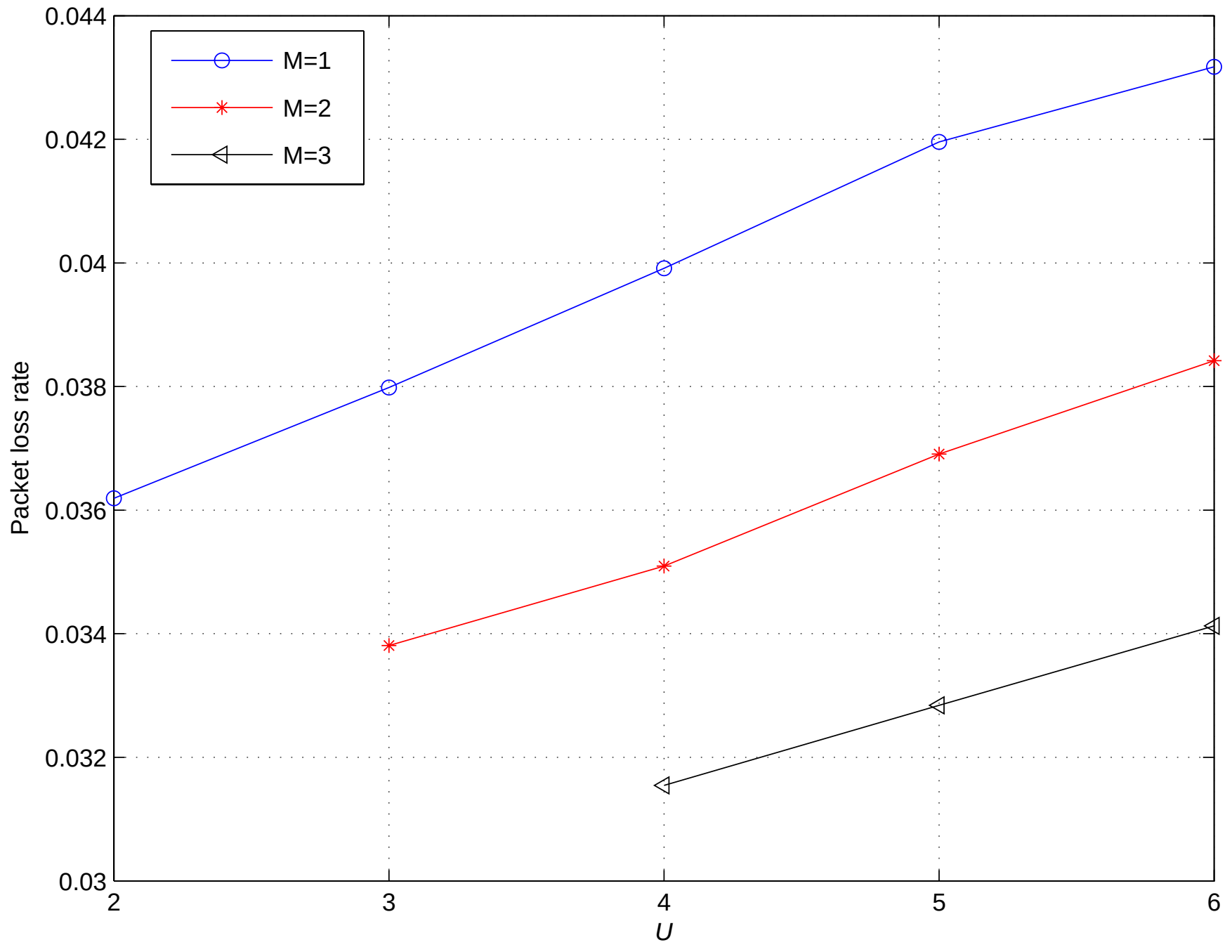

**Fig. 5:** Impact of *U* on packet loss rate

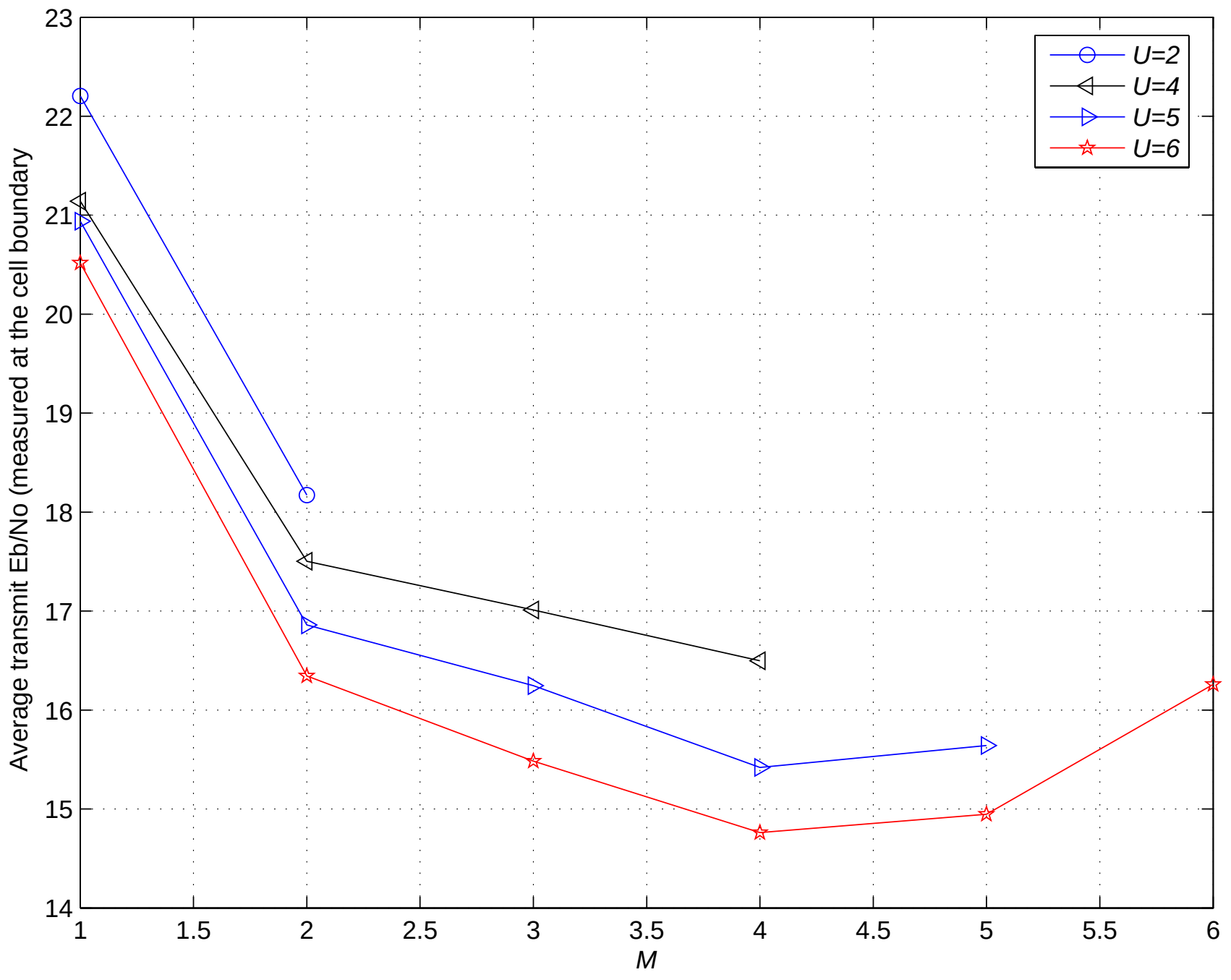

**Fig. 6:** Impact of *M* and *U* on power efficiency. *U*=6 and target BER is $10^{-6}$.

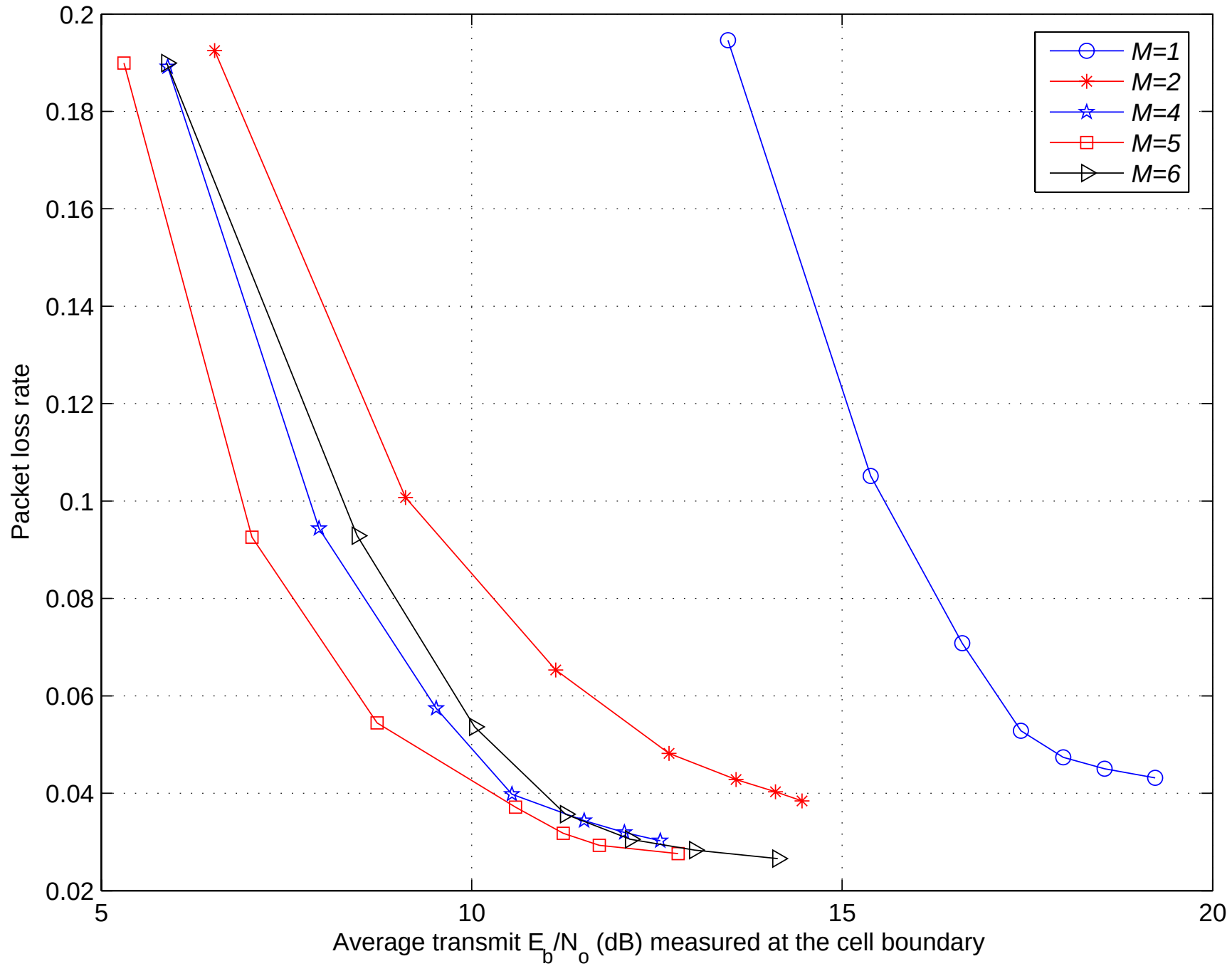

**Fig. 7:** Pareto optimum. For each *U*, there is an optimal *M*. In the figure, *U*=6.

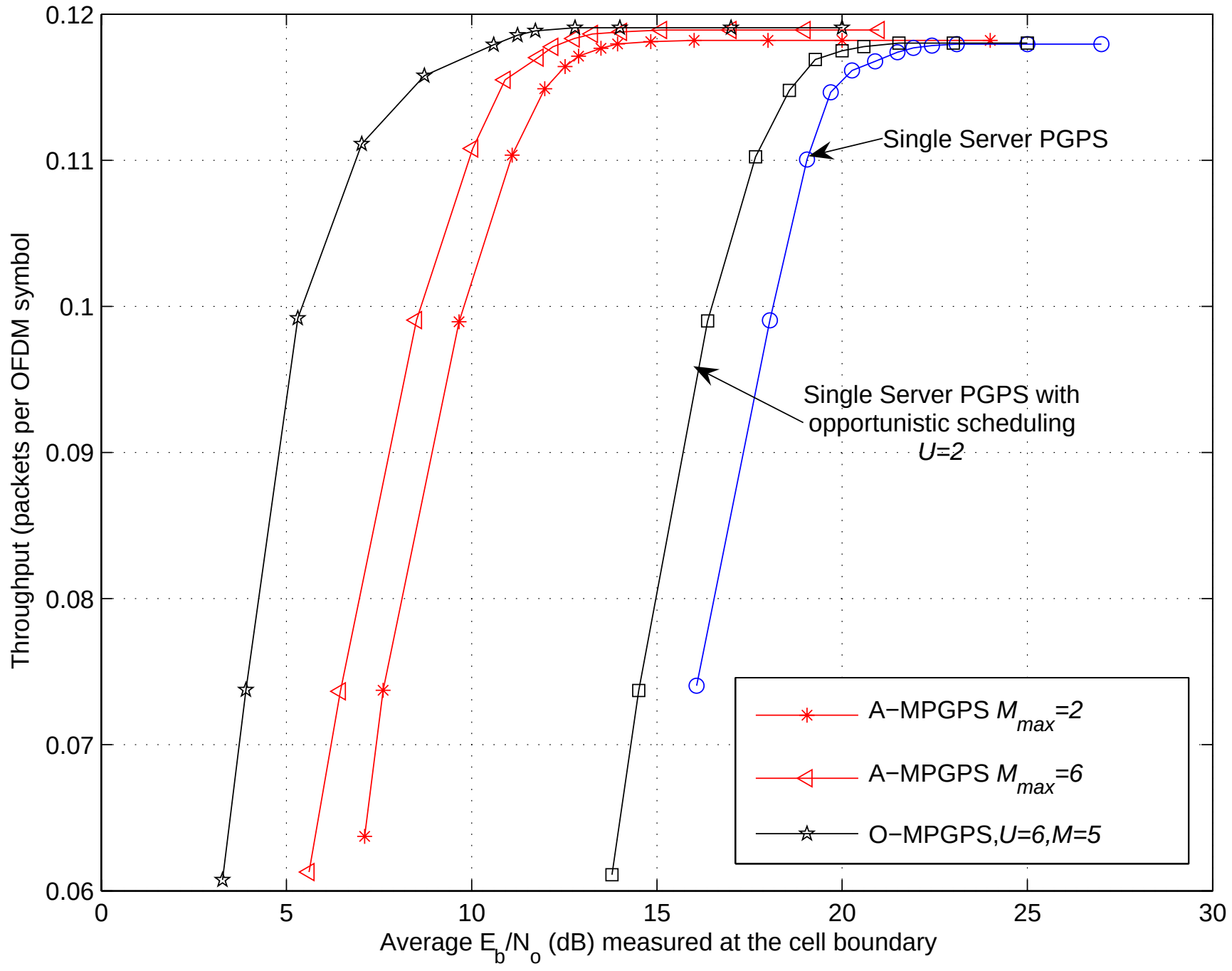

**Fig. 8:** Throughput comparison between A-MPGPS, O-MPGPS, and the conventional single-server PGPS schemes

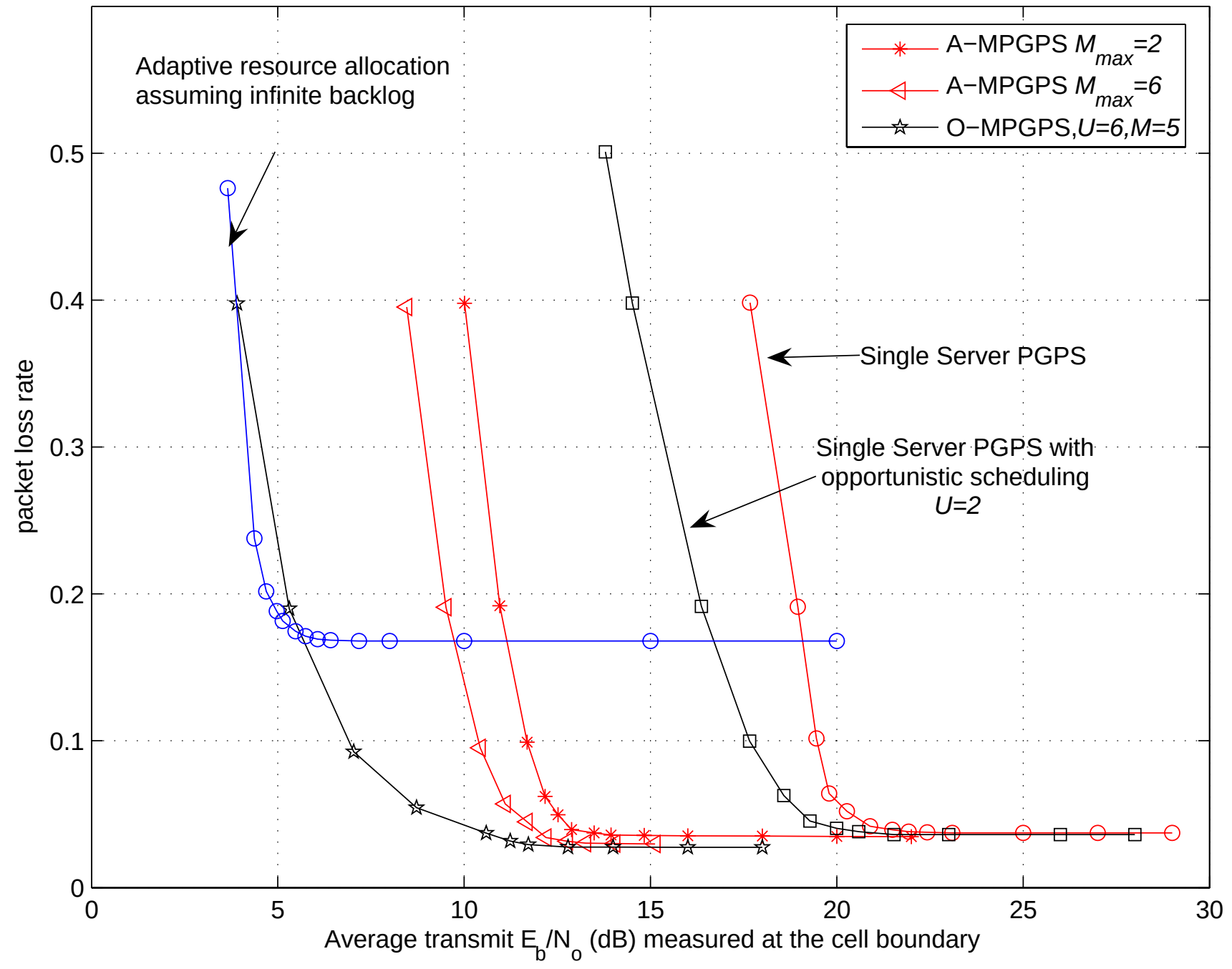

**Fig. 9:** Packet loss rate comparison between A-MPGPS, O-MPGPS, and the conventional single-server PGPS and resource allocation schemes